\documentclass[reqno,a4paper,12pt]{amsart} 

\usepackage{amsmath,amscd,amsfonts,amssymb}
\usepackage{mathrsfs,dsfont}

\usepackage{tikz}
\usetikzlibrary{arrows.meta,decorations.markings}
\definecolor{contourblue}{RGB}{31,90,166}
\tikzset{
	contour/.style={
		contourblue,
		very thick,
		line join=round,
		postaction={decorate},
		decoration={
			markings,
			mark=at position 0.55 with {\arrow{Latex[length=3mm,width=2mm]}}
		}
	},
	axis/.style={gray!70, thin, -{Latex[length=2.4mm,width=1.6mm]}},
	tick/.style={gray!70, thin}
}

\allowdisplaybreaks

\numberwithin{equation}{section}
\numberwithin{figure}{section}

\newcommand\C{\mathbb{C}}
\newcommand{\R}{\mathbb{R}}
\newcommand{\Z}{\mathbb{Z}}
\newcommand{\SSS}{\mathcal{S}}

\newcommand{\la}{\langle}
\newcommand{\ra}{\rangle}

\newcommand{\ft}[1]{\widehat{#1}}
\newcommand{\supp}{\operatorname{supp}}

\renewcommand\le{\leqslant}
\renewcommand\ge{\geqslant}

\theoremstyle{plain}
\newtheorem{thm}{Theorem}[section]

\newtheorem*{claim*}{Claim}

\theoremstyle{definition}

\newtheorem*{definition*}{Definition}
\newtheorem*{remarks*}{Remarks}
\newtheorem*{remark*}{Remark}

\newenvironment{enumerate-alph}
{\begin{enumerate}
		\addtolength{\itemsep}{5pt}
		}
	{\end{enumerate}}

\newenvironment{enumerate-num}
{\begin{enumerate}
		\addtolength{\itemsep}{5pt}
		}
	{\end{enumerate}}

\newenvironment{enumerate-text}
{\begin{enumerate}
		\addtolength{\itemsep}{5pt}
		}
	{\end{enumerate}}

\begin{document}
	
	\title
	[Crystalline-type measure and an application to translational tiling]
	{A new example of crystalline-type measure with unit masses and an application to tiling of the real line by translates of a function}
	
	\author{Anton Tselishchev}
	\address{St. Petersburg Department of Steklov Mathematical Institute, Fontanka 27, St. Petersburg 191023, Russia}
	\email{celis\_anton@pdmi.ras.ru}


	\begin{abstract}
		In this note we provide a simple explicit example of a discrete set $\Lambda\subset\mathbb{R}$ of unbounded density such that the Fourier transform of its counting measure $\widehat{\delta}_\Lambda$ is a tempered distribution of the form $c\delta_0'' + \mu$ where $\mu$ is a discrete measure supported by a union of three arithmetic progressions. This example is then used to construct a translational tiling of unbounded density of the real line  by a function with compact support.
	\end{abstract}

	\maketitle


	\section{Introduction}
	
	\subsection{A new crystalline-type measure with unit masses}
	
	A crystalline measure on $\R$ is a purely atomic measure $\mu$, supported by a discrete (i.e., locally finite) set $\Lambda$, which is a tempered distribution such that its distributional Fourier transform $\widehat{\mu}$ is also a measure supported by a discrete set $S$:
	\begin{equation}\label{eq:crystmes}
		\mu=\sum_{\lambda\in\Lambda} a_\lambda \delta_\lambda, \qquad \widehat{\mu} = \sum_{s\in S} b_s \delta_s.	
	\end{equation}
	This condition can be rewritten as a Fourier summation formula: for every Schwartz function $f$ on $\R$ 
	\begin{equation}\label{eq:summformula}
		\sum_{\lambda\in\Lambda} a_\lambda \widehat{f}(\lambda) = \sum_{s\in S} b_s f (s).
	\end{equation}
	
	The basic example of a crystalline measure is 
	\begin{equation}
		\mu = \delta_{\mathbb{Z}} = \sum_{n\in\mathbb{Z}}\delta_n.
	\end{equation}
	According to the Poisson summation formula we have $\widehat{\delta}_{\Z} = \delta_\Z$, therefore this measure satisfies the conditions $\eqref{eq:crystmes}$ with $\Lambda=S=\Z$. Here, and everywhere below, we use the following normalization for Fourier transform of a Schwartz function $f$: $\widehat{f}(t)=\int_\R e^{-2\pi i t x} f(x)\, dx$. As always, Fourier transform on the class of tempered distributions $\SSS'(\R)$ is defined by the equality $\la \widehat{u}, f \ra = \la u, \widehat{f} \ra$. Basic information about tempered distributions and their Fourier transform can be found e.g. in \cite[Chapter 2]{Gra14}.
	
	For the following discussion we need two definitions. A set $\Lambda\subset \R$ is called uniformly discrete if there exists $\delta > 0$ such that the distance between any two distinct points of $\Lambda$ is at least $\delta$. Next, a set $\Lambda$ has bounded density if it is a finite union of uniformly discrete sets. Equivalently, $\Lambda$ has bounded density if there exists $C>0$ such that
	\begin{equation}
		\# (\Lambda\cap [x, x+1)) \le C
	\end{equation}
	for every $x\in\R$.
	
	Most recent results on the topic of crystalline measures fall into one of two categories: they either prove that under certain additional assumptions on the sets $S$ and $\Lambda$ any formula \eqref{eq:summformula} can be obtained from the Poisson summation formula by trivial operations or they provide new non-trivial examples of such formulas. In the first direction an important result was obtained in \cite{LO15}: it was shown that if both sets $\Lambda$ and $S$ are uniformly discrete then $\Lambda$ must be contained in a finite union of translates of a lattice, the same is true for $S$ and the formula \eqref{eq:summformula} can be obtained from the Poisson summation formula in a certain trivial way (we do not present the exact rigorous formulation of the results from \cite{LO15} here); these results were later generalized to the setting of crystalline tempered distributions in \cite{LR21}. In the other direction, constructions of non-trivial crystalline measures with unit masses $a_\lambda$ whose support $\Lambda$ is uniformly discrete but is not contained in any finite union of arithmetic progressions were obtained, by using different approaches, in \cite{KS20}, \cite{Mey23}, \cite{OU20}. This short list of results is, of course, very much incomplete, we do not survey all progress in this topic here.
	
	For a discrete set $\Lambda\subset \R$ we denote by $\delta_\Lambda$ the measure with unit masses at the points of $\Lambda$, i.e., $\delta_\Lambda = \sum_{\lambda\in\Lambda}\delta_\lambda$.
	
	Our first result, motivated by an application to an open problem about tilings of the real line by translates of a function (we discuss this application below) provides a new crystalline-type measure $\delta_\Lambda$. By ``crystalline-type measure'' we mean here that $\widehat{\delta_\Lambda}$ is not exactly of the form \eqref{eq:crystmes} because it contains also the term $\delta_0''$. If we adopt the terminology from \cite{LR21}, we can say that $\delta_\Lambda$ is a measure which is simultaneously a crystalline tempered distribution.
	
	We denote by $\Z_+$ the set of nonnegative integer numbers. Put 
	\begin{equation}
		\omega_1 = 1, \ \omega_2 = \sqrt{2},  \ \omega_3 = \sqrt{3}, \quad a=\frac{\omega_1 + \omega_2 + \omega_3}{2}
	\end{equation}
	and consider the sets
	\begin{equation}\label{eq:lambdadef}
		\Lambda_1 = \{a+ n_1\omega_1 + n_2\omega_2 + n_3\omega_3: n_1, n_2, n_3\in\Z_+\}, \qquad \Lambda = \Lambda_1 \cup (-\Lambda_1).
	\end{equation}
	Denote also 
	\begin{equation}\label{eq:P_and_Q}
		P = \omega_1\omega_2\omega_3 = \sqrt{6},\qquad Q = \omega_1^2 + \omega_2^2 + \omega_3^2 = 6.
	\end{equation} Notice that our set $\Lambda$ is not uniformly discrete and it does not have bounded density but it has temperate growth: we have
	\begin{equation}\label{eq:lambdagrowth}
		\# (\Lambda\cap [-R, R]) = O(R^3), \qquad R\to\infty.
	\end{equation}
	Therefore, $\delta_\Lambda$ is a tempered distribution.
	
	\begin{thm} \label{thm:cryst}
		The following identity holds in $\SSS'(\R)$\emph{:}
		\begin{equation}\label{eq:mainthm}
			\widehat{\delta}_\Lambda = -\frac{1}{8\pi^2 P}\delta_0'' - \frac{Q}{24P} \delta_0 + \sum_{k\neq 0} (b_k\delta_k + c_k\delta_{k/\sqrt{2}} + d_k \delta_{k/{\sqrt{3}}}),
		\end{equation}
		where coefficients $b_k$, $c_k$, $d_k$ are given by the formulas
		\begin{align}
			\label{eq:b_k} b_k &= \frac{(-1)^{k+1}}{4\omega_1 \sin(\pi k \frac{\omega_2}{\omega_1}) \sin (\pi k  \frac{\omega_3}{\omega_1})}, \\
			\label{eq:c_k} c_k &= \frac{(-1)^{k+1}}{4\omega_2\sin (\pi k  \frac{\omega_1}{\omega_2}) \sin(\pi k  \frac{\omega_3}{\omega_2})},\\
			\label{eq:d_k} d_k &= \frac{(-1)^{k+1}}{4\omega_3 \sin(\pi k  \frac{\omega_1}{\omega_3}) \sin (\pi k  \frac{\omega_2}{\omega_3})}.
		\end{align}
		The series in \eqref{eq:mainthm} converges unconditionally in $\SSS'(\R)$ (i.e., after testing against any Schwartz function we get an absolutely convergent series).
		
		In particular,
		\begin{equation}\label{eq:Sdef}
			S:= \mathrm{supp}\, \ft{\delta}_\Lambda=\Z \cup \Big(\frac{1}{\sqrt{2}}\Z\Big) \cup \Big(\frac{1}{\sqrt{3}}\Z\Big).
		\end{equation}
	\end{thm}
	
	\textbf{Remarks.} 1. Of course, one can take other values of $\omega_1, \omega_2, \omega_3$, however we do not claim that for \emph{any} choice of rationally independent numbers the corresponding series in \eqref{eq:mainthm} converges in $\SSS'(\R)$: for certain choices of these numbers the coefficients $b_k$, $c_k$, $d_k$ may have rapid growth. In this case the corresponding Fourier summation formula should be formulated more accurately, for a narrower class of functions than $\SSS(\R)$. We do not address this issue here. 
	
	2. It is also possible to obtain similar formulas (with higher-order derivatives on the right-hand side of \eqref{eq:mainthm}) by taking odd number of rationally independent numbers instead of three numbers $\omega_1$, $\omega_2$, $\omega_3$. A certain formula of this kind can also be obtained by taking even number of rationally independent numbers but in this case instead of $\delta_\Lambda$ one would have to take a measure with charges $\pm 1$ at the points of the corresponding set $\Lambda$.
	
	3. Our proof below uses the contour integration approach from \cite{OU20}. Notice that our set $\Lambda$ would be $S$ in notation of the paper \cite{OU20}, and vice versa. In this sense our example of ``crystalline-type measure'' $\delta_\Lambda$ is not completely similar to examples obtained in \cite{OU20}: in notation of \cite[Corollary 1]{OU20} we manage to apply approach from that paper so that all masses of $\ft{\mu}$, not $\mu$, are unit.
	
	4. Theorem~\ref{thm:cryst} gives another example of the set $\Lambda$ which satisfies the conclusion of Theorem 5.1 in \cite{KL21}. In that theorem the authors' approach allowed only the control of $\ft{\delta}_\Lambda$ on the interval $(-\frac{1}{2}, \frac{1}{2})$. Our construction is completely different and it allows to precisely compute $\ft{\delta}_{\Lambda}$.
	
	\subsection{An application to translational tiling by a function}
	
	Suppose that $f\in L^1(\R)$ and $\Lambda\subset\R$ is a discrete set. The function $f$ tiles $\R$ at level $w\in\mathbb{C}$ with translation set $\Lambda$ if
	\begin{equation}
		\sum_{\lambda\in\Lambda} f(x-\lambda) = w \quad \text{a.e.}
	\end{equation}
	and the series in the above formula converges absolutely a.e. We refer the reader to the paper \cite{KL21} for a survey of the results about tilings by translates of a function.
	
	One of the main structural theorems, obtained in \cite{LM91} and then independently in \cite{KL96}, states that if a nonzero function $f\in L^1(\R)$ has compact support and tiles $\R$ at some level $w$ with translation set $\Lambda$ of bounded density then such tiling is necessarily periodic, i.e. $\Lambda$ is a finite union of disjoint arithmetic progressions.
	
	It is known that the assumption on the support of $f$ is essential in this result since in general non-periodic tilings by translates of a function do exist, an example of such tiling was constructed in \cite{KL16}. We note, however, that it is not known whether there exist examples of non-periodic tilings where support of the function $f$ has finite measure.
	
	Not much is known about tilings when translation set $\Lambda$ does not have bounded density. The first example of such tiling by translates of a Schwartz function was constructed in \cite{KL21}. The Schwartz function $f$ in this construction, however, has an unbounded support. It has been therefore an open problem, formulated in \cite[Section 6.3]{KL21}, whether there exist tilings by functions with compact support with translation set of unbounded density. In other words, this problem asks whether the ``bounded density'' assumption on $\Lambda$ is crucial in the aforementioned result from \cite{LM91} and \cite{KL96}.
	
	A simple application of Theorem~\ref{thm:cryst} provides an answer to this open problem: such translational tilings by a function with compact support with the translation set of unbounded density do exist.
	
	\begin{thm}\label{thm:tiling}
		Suppose that the set $\Lambda$ is given by the formula \eqref{eq:lambdadef} and $w$ is a complex number. Then there exists a non-zero function $f\in\SSS (\R)$ with compact support such that
		\begin{equation}\label{eq:tiling_ident}
			\sum_{\lambda\in\Lambda} f(x-\lambda) = w \quad \text{for every}\ x\in\R.
		\end{equation}
	\end{thm}
	Notice that since the set $\Lambda$ is discrete and $f$ has compact support, the series in the formula \eqref{eq:tiling_ident} is automatically absolutely convergent for every $x\in\R$: in fact, it has only finite number of nonzero summands.
	
	\textbf{Remark.} Our Fourier-analytic proof of Theorem~\ref{thm:tiling} can be reformulated in elementary terms without using Fourier transform: in fact, an example of a function $f\in\SSS(\R)$ which tiles $\R$ with the set $\Lambda$ can be obtained from an arbitrary Schwartz function $g$ with compact support as a certain (finite) linear combination of its translates; after writing down the formulas, the tiling condition can be verified directly. However, an approach which uses Fourier analysis is probably more natural in this problem and Theorem~\ref{thm:cryst} seems to be interesting in its own right.
	
	In the next two sections we prove Theorems~\ref{thm:cryst} and \ref{thm:tiling}
	
	\subsection{Usage of Large Language Models}
	
	An example from Theorem~\ref{thm:cryst} was found by ChatGPT 5.6 Sol (Pro). The model's proof, however, was different from the one presented in the paper below; the proof in this paper, which uses an approach to the construction of crystalline measures from \cite{OU20}, belongs to the author. The application to the tiling problem, which was the initial goal and motivation for searching for such kind of formulas, is also due to the author.
	
	The TikZ code for Figure~\ref{fig_contour} below was generated by ChatGPT. Apart from that, the text of the paper was written manually by the author without using LLMs (therefore at least all possible typos are the author's contribution). The author takes full responsibility for the mathematical content of the paper.
	
	\section{Computation of $\ft{\delta}_\Lambda$: the proof of Theorem~\ref{thm:cryst}}
	
	\subsection{}
	
	At first we show that the right-hand side of the formula \eqref{eq:mainthm} defines a tempered distribution and the series in this formula is unconditionally convergent in $\SSS'(\R)$. It is enough to show that for any function $f\in \SSS(\R)$ all three series
	\begin{equation}
		\sum_{k\neq 0} b_k f(k), \quad \sum_{k\neq 0} c_k f(k/\sqrt{2}), \quad \sum_{k\neq 0} d_k f(k/\sqrt{3})
	\end{equation}
	are absolutely convergent and define continuous linear functionals on the space $\SSS (\R)$. This in turn would follow from the estimates
	\begin{equation}
		|b_k| \lesssim k^2, \quad |c_k|\lesssim k^2, \quad |d_k|\lesssim k^2.
	\end{equation}
	We prove the first of these estimates since the proof of the rest are similar.
	
	Recall that for our concrete numbers $\omega_1, \omega_2, \omega_3$ the coefficients $b_k$ are defined as
	\begin{equation}
		b_k = \frac{(-1)^{k+1}}{4\sin(\pi k \sqrt{2})\sin(\pi k \sqrt{3})}.
	\end{equation}
	Since $|\sin(\pi x)|\asymp \mathrm{dist}\, (x, \Z)$ for $x\in\R$, we only need to show that
	\begin{equation}
		\mathrm{dist}\, (k\sqrt{2}, \Z) \gtrsim \frac{1}{k} \quad \text{and}\quad  \mathrm{dist}\, (k\sqrt{3}, \Z) \gtrsim \frac{1}{k}.
	\end{equation}
	Again, these estimates can be proved in a similar way, so we show the proof only of the first one. It is equivalent to the fact that for every integer $m$ we have
	\begin{equation}
		|k\sqrt{2} - m|\gtrsim \frac{1}{k}.
	\end{equation}
	Without loss of generality we can assume that $k > 0$ and then we can also assume that $m \ge 0$ and $|m-k\sqrt{2}|\le 1$ and therefore $0\le m\le 3k$ (otherwise the estimate is trivial). The left-hand side of the above inequality equals
	\begin{equation}
		\Big| \frac{2k^2-m^2}{m+k\sqrt{2}} \Big|.
	\end{equation}
	Obviously, $m+k\sqrt{2}\le 5k$ and $|2k^2-m^2| \ge 1$ because it is a nonnegative integer and it cannot be equal to zero by irrationality of $\sqrt{2}$. These simple observations prove the required estimate.
	
	\subsection{}
	
	Compactly supported functions are dense in $\SSS(\R)$, see e.g. \cite[Theorem 7.10]{Rud91}. Since Fourier transform is an automorphism of $\SSS(\R)$, the same is true for functions with compactly supported Fourier transform. Therefore, if we denote by $u$ the distribution on the right-hand side of \eqref{eq:mainthm}, we need to show that 
	\begin{equation}
		\langle \ft{\delta}_\Lambda , F \rangle = \langle u, F \rangle	
	\end{equation}
	for any function $F\in \SSS(\R)$ such that $f=\ft{F}$ has compact support. In other words, we need to show that the following identity for any such function $F$:
	\begin{equation}\label{eq:toprove_main}
		\sum_{\lambda\in\Lambda} \ft{F}(\lambda) = \frac{-F''(0)}{8\pi^2 P} - \frac{Q F(0)}{24 P} + \sum_{k\neq 0} (b_k F(k) + c_k F(k/\sqrt{2}) + d_k F(k/\sqrt{3})).
	\end{equation}
	We fix a function $F=f^\vee$ such that $f$ is a Schwartz function with compact support and our goal is to prove the above formula for it.
	
	Suppose that $\supp f\subset [-R, R]$. The function $F$ is an inverse Fourier transform of $f$ and can be extended as an entire function to the whole complex plane $\C$ in a standard way:
	\begin{equation}\label{eq:fourier_transform}
		F(z) = \int_\R e^{2\pi i x z} f(x)\, dx.
	\end{equation}
	We note that $F$ belongs to the Paley--Wiener space $\mathrm{PW}_R$; see e.g. \cite[Lecture 2]{OU16} for more information about Paley--Wiener spaces. In particular, the following estimate with constant $C_1$ depending on $f$ can be easily verified by putting absolute value inside the integral in \eqref{eq:fourier_transform}
	\begin{equation}
		|F(z)| \le C_1 e^{2\pi R |\mathrm{Im}\, z|}.
	\end{equation}
	Besides that, integrating by parts we get
	\begin{equation}
		-4\pi^2z^2 F(z) = \int_\R e^{2\pi i x z} f''(x)\, dx,
	\end{equation}
	and hence, similarly to the above we infer that
	\begin{equation}\label{eq:F_decay}
		|F(z)|\le C_2 \frac{1}{|z|^2} e^{2\pi R |\mathrm{Im}\, z|}.
	\end{equation}
	In particular, for any fixed horizontal strip in the complex plane there exists a constant $C_3 > 0$ such that for any $z$ in this strip we have
	\begin{equation}\label{eq:F_decay_strip}
		|F(z)|\le\frac{C_3}{|z|^2}.
	\end{equation}
	
	Consider the following function:
	\begin{equation}
		\Phi(z) = \sin (\pi\omega_1 z)\cdot \sin(\pi\omega_2 z)\cdot \sin(\pi\omega_3 z).
	\end{equation}
	Then
	\begin{align}
		\label{eq:Phi_1}	\frac{1}{\Phi(z)} &= 8i
		\frac{e^{2\pi i a z}}{(1-e^{2\pi i \omega_1 z})(1-e^{2\pi i \omega_2 z})(1-e^{2\pi i \omega_3 z})}
		\\ \label{eq:Phi_2} &= -8i  \frac{e^{-2\pi i a z}}{(1-e^{-2\pi i \omega_1 z})(1-e^{-2\pi i \omega_2 z})(1-e^{-2\pi i \omega_3 z})}.
	\end{align}
	
	Now we choose positive real numbers $M_n$ tending to infinity so that  \begin{equation}\label{eq:Phi_est}
		|\Phi(z)|\ge\delta\qquad \text{for all}\ z\ \text{such that} \ \mathrm{Re}\, z = \pm M_n
	\end{equation} 
	with some constant $\delta > 0$. It is easy to see that such choice exists since the entire function $\sin(z)$ is uniformly bounded away from zero when $z$ is bounded away from the points $\{\pi k\}_{k\in\Z}$: indeed, it is enough to ensure that all points $M_n$ are, say, at least 1/100 of the distance to the set $S$ given by the formula \eqref{eq:Sdef}.
	
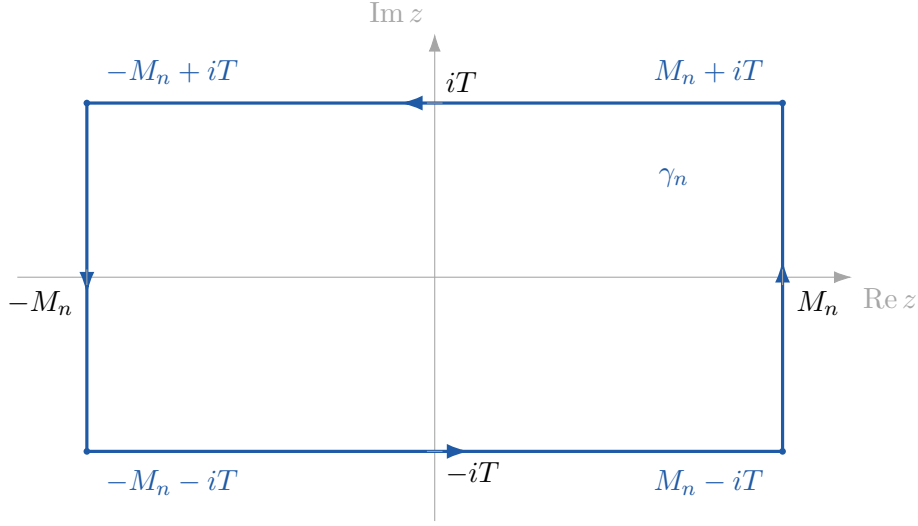
\begin{figure}
		\centering
		\begin{tikzpicture}[x=1.15cm,y=1.15cm,font=\small]
			
			\draw[axis] (-4.8,0) -- (4.8,0) node[below right] {$\operatorname{Re} z$};
			\draw[axis] (0,-2.8) -- (0,2.8) node[above left] {$\operatorname{Im} z$};

			\draw[contour] (-4,-2) -- (4,-2);
			\draw[contour] (4,-2) -- (4,2);
			\draw[contour] (4,2) -- (-4,2);
			\draw[contour] (-4,2) -- (-4,-2);

			\fill[contourblue] (-4,-2) circle (1.25pt);
			\fill[contourblue] (4,-2) circle (1.25pt);
			\fill[contourblue] (4,2) circle (1.25pt);
			\fill[contourblue] (-4,2) circle (1.25pt);

			\draw[tick] (-4,0.09) -- (-4,-0.09);
			\node[below left=4pt,fill=white,inner sep=1pt] at (-4,0) {$-M_n$};
			\draw[tick] (4,0.09) -- (4,-0.09);
			\node[below right=4pt,fill=white,inner sep=1pt] at (4,0) {$M_n$};
			\draw[tick] (-0.09,2) -- (0.09,2);
			\node[above right=3pt,fill=white,inner sep=1pt] at (0,2) {$iT$};
			\draw[tick] (-0.09,-2) -- (0.09,-2);
			\node[below right=3pt,fill=white,inner sep=1pt] at (0,-2) {$-iT$};
			
			\node[above right=3pt,text=contourblue] at (-4,2)
			{$-M_n+iT$};
			\node[above left=3pt,text=contourblue] at (4,2)
			{$M_n+iT$};
			\node[below right=3pt,text=contourblue] at (-4,-2)
			{$-M_n-iT$};
			\node[below left=3pt,text=contourblue] at (4,-2)
			{$M_n-iT$};
			
			\node[contourblue,fill=white,inner sep=2pt] at (2.75,1.15)
			{$\gamma_n$};
		\end{tikzpicture}
		\caption{Contour $\gamma_n$.}
		\label{fig_contour}
	\end{figure}

	Fix an arbitrary number $T > 0$ and consider the contour $\gamma_n$ which is a rectangle with vertices $\pm M_n\pm iT$ oriented counterclockwise in $\C$ (see Figure~\ref{fig_contour}). The main idea of the proof is to compute the quantity
	\begin{equation}
		\lim_{n\to\infty} \int_{\gamma_n} \frac{F(z)}{\Phi(z)}\, dz
	\end{equation}
	in two different ways: directly and using the Residue Theorem.
	
	\subsection{Direct computation of the integral}
	
	Using \eqref{eq:F_decay_strip} and \eqref{eq:Phi_est}, it is easy to conclude that the integrals over two lateral sides of the rectangle $\gamma_n$ (that is, over the segments $[\pm M_n - iT, \pm M_n + iT]$) tend to zero. Hence we need to compute
	\begin{equation}\label{eq:2.19}
		\int_{-\infty - iT}^{\infty - iT}  \frac{F(z)}{\Phi(z)}\, dz - \int_{-\infty + iT}^{\infty + iT}  \frac{F(z)}{\Phi(z)}\, dz.
	\end{equation}
	Here both integrals converge absolutely: it follows from the decay of the function $F$ given by \eqref{eq:F_decay_strip} and the fact that $\Phi(z)$ is uniformly bounded from zero for all $z$ such that $\mathrm{Im}\, z = \pm T$.
	
	We first show how to compute the second integral in \eqref{eq:2.19}. We use formula \eqref{eq:Phi_1} for the function $\frac{1}{\Phi(z)}$ and expand it using the geometric series:
	\begin{equation}
		\frac{1}{\Phi(z)} = 8i e^{2\pi i a z} \cdot \Big( \sum_{n_1 = 0}^\infty e^{2\pi i n_1\omega_1 z} \Big)\cdot \Big( \sum_{n_2 = 0}^\infty e^{2\pi i n_2\omega_2 z} \Big)\cdot \Big( \sum_{n_3 = 0}^\infty e^{2\pi i n_3\omega_3 z} \Big).
	\end{equation}
	All series here absolutely uniformly converge for $z$ such that $\mathrm{Im}\, z = T$. Hence, substituting this formula into our integral and opening the brackets, we get, according to the definition of the set $\Lambda_1$ given by \eqref{eq:lambdadef}, that 
	\begin{equation}\label{eq:2.21}
		\int_{-\infty + iT}^{\infty + iT}  \frac{F(z)}{\Phi(z)}\, dz = 8i \sum_{\lambda\in \Lambda_{1}} \int_{\R + iT} F(z) e^{2\pi i \lambda z}\, dz.
	\end{equation}
	Finally, using a simple contour integration argument one can show that for any $\lambda\ge 0$
	\begin{equation}
		\int_{\R + iT} F(z) e^{2\pi i \lambda z}\, dz = \int_{\R} F(z) e^{2\pi i \lambda z}\, dz=\ft{F}(-\lambda).
	\end{equation}
	(In order to prove it, one can notice that the integral over the rectangle with vertices $\pm L$ and $\pm L + iT$ of the holomorphic function $ F(z) e^{2\pi i \lambda z}$ is equal to zero and then let $L\to\infty$ and again use the decay condition \eqref{eq:F_decay_strip} and the fact that $|e^{2\pi i \lambda z}|\le 1$ on this contour; we omit the routine details.)
	
	Since $\ft{F}$ has compact support, only a finite number of summands in the right-hand side of \eqref{eq:2.21} is nonzero and we have shown that
	\begin{equation}
		\int_{-\infty + iT}^{\infty + iT}  \frac{F(z)}{\Phi(z)}\, dz = 8i \sum_{\lambda\in (-\Lambda_1)} \ft{F}(\lambda).
	\end{equation}
	
	The first integral in \eqref{eq:2.19} can be computed in a similar way: we only have to use the formula \eqref{eq:Phi_2} instead of \eqref{eq:Phi_1} for the function $\frac{1}{\Phi (z)}$ in order to expand it as the product of geometric series convergent in the lower half-plane. Then similar computations then yield
	\begin{equation}
		\int_{-\infty - iT}^{\infty - iT}  \frac{F(z)}{\Phi(z)}\, dz = -8i  \sum_{\lambda\in\Lambda_1} \ft{F} (\lambda).
	\end{equation}
	
	Summing up, we have shown that 
	\begin{equation}\label{eq:integral_direct}
		\lim_{n\to\infty} \int_{\gamma_n} \frac{F(z)}{\Phi(z)}\, dz = -8i \sum_{\lambda\in\Lambda} \ft{F} (\lambda).
	\end{equation}

	\subsection{Computation of the integral using the Residue Theorem}
	
	The function $\frac{F(z)}{\Phi(z)}$ is meromorphic in $\mathbb{C}$ and its poles are located in the set
	\begin{equation}
		S=\Big(\frac{1}{\omega_1}\Z\Big) \cup \Big(\frac{1}{\omega_2}\Z\Big) \cup \Big(\frac{1}{\omega_3}\Z\Big) = \Z \cup \Big(\frac{1}{\sqrt{2}}\Z\Big) \cup \Big(\frac{1}{\sqrt{3}}\Z\Big).
	\end{equation}
	Hence by the Residue Theorem we have
	\begin{equation}\label{eq:resthm}
		\lim_{n\to\infty} \int_{\gamma_n} \frac{F(z)}{\Phi(z)}\, dz =2\pi i \sum_{w\in S} \mathrm{Res}\Big(\frac{F}{\Phi}, w \Big).
	\end{equation}
	The poles of the function $\frac{1}{\Phi}$ at all points of the form $k/\omega_j$ with $k\neq 0$
	are of the first order and the required residues can be computed in a very simple way. For example,
	\begin{equation}
		\mathrm{Res}\Big(\frac{F}{\Phi}, \frac{k}{\omega_1} \Big) = \frac{F(\frac{k}{\omega_1})}{\sin (\pi k \frac{\omega_2}{\omega_1})\sin (\pi k \frac{\omega_3}{\omega_1})}\cdot \frac{1}{\pi \omega_1 \cos(\pi k)}=- \frac{4}{\pi} b_k F(k),
	\end{equation}
	where $b_k$ is given by \eqref{eq:b_k}.
	
	Similarly, 
	\begin{equation}
		\mathrm{Res}\Big(\frac{F}{\Phi}, \frac{k}{\omega_2} \Big) = -\frac{4}{\pi} c_k F\Big(\frac{k}{\sqrt{2}}\Big), \qquad  \mathrm{Res}\Big(\frac{F}{\Phi}, \frac{k}{\omega_3} \Big) = -\frac{4}{\pi} d_k F\Big(\frac{k}{\sqrt{3}}\Big),
	\end{equation}
	where the coefficients $c_k$ and $d_k$ are given by \eqref{eq:c_k} and \eqref{eq:d_k}.
	
	It remains only to compute the residue of our function at zero. The probably easiest way to quickly do it is to use Laurent series
	\begin{equation}
		\frac{1}{\sin z} = \frac{1}{z} + \frac{1}{6} z + \ldots
	\end{equation}
	as follows:
	\begin{align}
		\frac{F(z)}{\Phi(z)} = \Big(F(0) + F'(0)z + \frac{1}{2} F''(0)z^2 + \ldots\Big)\\
		\times  \Big(\frac{1}{\pi \omega_1 z} + \frac{1}{6}\pi\omega_1 z + \ldots\Big)\Big(\frac{1}{\pi \omega_2 z} + \frac{1}{6}\pi\omega_2 z + \ldots\Big)\Big(\frac{1}{\pi \omega_3 z} + \frac{1}{6}\pi\omega_3 z + \ldots\Big).
	\end{align}
	Opening the brackets and considering the coefficient of $z^{-1}$ in the resulting series gives us
	\begin{equation}
		\mathrm{Res}\Big(\frac{F}{\Phi}, 0 \Big) = \frac{F''(0)}{2\pi^3\omega_1\omega_2\omega_3} + \frac{F(0)}{6\pi}\cdot \Big( \frac{\omega_1}{\omega_2\omega_3} + \frac{\omega_2}{\omega_1\omega_3} + \frac{\omega_3}{\omega_1\omega_2} \Big)= \frac{F''(0)}{2\pi^3 P} + \frac{QF(0)}{6\pi P},
	\end{equation}
	where, as in \eqref{eq:P_and_Q}, $P = \omega_1\omega_2\omega_3 = \sqrt{6}$ and  $Q = \omega_1^2 + \omega_2^2 + \omega_3^2 = 6$.
	
	Substituting these computations into \eqref{eq:resthm}, we get that
	\begin{align}
		&\lim_{n\to\infty} \int_{\gamma_n} \frac{F(z)}{\Phi(z)}\, dz \\
		= i\frac{F''(0)}{\pi^2 P} + i\frac{QF(0)}{3 P} -8i &\Big(  \sum_{k\neq 0} (b_k F(k) + c_k F(k/\sqrt{2}) + d_k F(k/\sqrt{3}))\Big).
	\end{align}
	It remains to compare this formula with \eqref{eq:integral_direct} and to conclude that the required identity \eqref{eq:toprove_main} is proved.
	
	\section{Application to translational tiling problem: the proof of Theorem~\ref{thm:tiling}}
	
	We now show how to apply the construction obtained in Theorem~\ref{thm:cryst} to construct an example of a tiling of unbounded density by a function with compact support. That is, we prove Theorem~\ref{thm:tiling}.
	
	Obviously, since we can always multiply our function $f$ by a nonzero constant, there exist essentially two different cases: $w=0$ and $w\neq 0$.
	
	\subsection{Case $w = 0$}
	
	Take an arbitrary nonzero Schwartz function $g$ with compact support and put 
	\begin{equation}\label{eq:hdef}
		h(t) = \ft{g}(t)\cdot  \sin (\pi\omega_1 t)\cdot \sin(\pi\omega_t z)\cdot \sin(\pi\omega_3 t).
	\end{equation}
	Clearly, $h$ is also a nonzero Schwartz function. Applying inverse Fourier transform to both sides of the above formula, since $(\sin t)^\vee$ is a combination of two delta-functions, one can conclude that $f=h^\vee$ is a linear combination of translates of $g$ and therefore it has compact support. The function $h$ also vanishes at the points of the set $S$ defined by \eqref{eq:Sdef} and it has zero of order at least 3 at the origin.
	
	Therefore, the formula \eqref{eq:mainthm} (or a direct application of the formula \eqref{eq:toprove_main}) implies that for $h = \ft{f}$ we have
	\begin{equation}
		h\cdot \ft{\delta}_{\Lambda} = 0.
	\end{equation}
	We can now take inverse Fourier transform and get that
	\begin{equation}
		f\ast \delta_\Lambda = 0,
	\end{equation}
	which can be equivalently rewritten as
	\begin{equation}
		\sum_{\lambda\in\Lambda} f(x-\lambda) = 0 \quad \text{for every}\ x\in\R.
	\end{equation}
	This finishes the proof for the case $w = 0$.
	
	\subsection{Case $w\neq 0$}
	
	This case can be treated similarly to the previous one. We can start with an arbitrary nonzero Schwartz function $g$ with compact support and assume additionally that $\ft{g}(0)\neq 0$. Then the function $h$ given by \eqref{eq:hdef} has zero of order exactly 3 at the origin. 
	
	The function $f$ can be defined as before by $f=h^\vee$ and then we put 
	\begin{equation}
		f_1(x) = \int_{-\infty}^x f(s)\, ds.
	\end{equation}
	Since the function $f\in\SSS(\R)$ compactly supported and has zero integral, the function $f_1$ is also a compactly supported Schwartz function. Moreover,
	\begin{equation}
		f_1'(x)  = f(x),
	\end{equation}
	hence
	\begin{equation}
		\ft{f}_1 (t) = \frac{\ft{f} (t)}{2\pi i t}.
	\end{equation}
	It means that $\ft{f}_1$ has zero of order exactly 2 at the origin and still vanishes at all points of the set $S$. Then the formula \eqref{eq:mainthm} implies that
	\begin{equation}
		\ft{f}_1\cdot \ft{\delta}_{\Lambda} = w\delta_0.
	\end{equation}
	for a nonzero number $w$. As before, we can now take inverse Fourier transform and get that
	\begin{equation}
		f_1\ast \delta_\Lambda = w,
	\end{equation}
	which can be equivalently rewritten as
	\begin{equation}
		\sum_{\lambda\in\Lambda} f_1(x-\lambda) = w \quad \text{for every}\ x\in\R,
	\end{equation}
	and we are done.
	

\end{document}